\documentclass{amsart}

\usepackage[utf8]{inputenc} 
\usepackage[english]{babel} 
\usepackage{mathrsfs}

\theoremstyle{plain} 
\newtheorem{thm}{Theorem} 
\newtheorem{lem}[thm]{Lemma}

\theoremstyle{definition}

\theoremstyle{remark} 

\newtheorem{que}{Question}

\begin{document} 
\title[Nehari's Theorem for Multiplicative Hankel Forms in $\mathcal{S}_p$]{Failure of Nehari's Theorem for Multiplicative Hankel Forms in Schatten Classes}
\date{\today} 

\author{Ole Fredrik Brevig} 
\address{Department of Mathematical Sciences, Norwegian University of Science and Technology (NTNU), NO-7491 Trondheim, Norway} 
\email{ole.brevig@math.ntnu.no}

\author{Karl-Mikael Perfekt}
\address{Department of Mathematical Sciences, Norwegian University of Science and Technology (NTNU), NO-7491 Trondheim, Norway}
\email{karl-mikael.perfekt@math.ntnu.no}

\thanks{The first author is supported by Grant 227768 of the Research Council of Norway.}



% ABSTRACT
\begin{abstract}
	Ortega-Cerd{\`a}--Seip demonstrated that there are bounded multiplicative Hankel forms which do not arise from bounded symbols. On the other hand, when such a form is in the Hilbert--Schmidt class $\mathcal{S}_2$, Helson showed that it has a bounded symbol. The present work investigates forms belonging to the Schatten classes between these two cases. It is shown that for every $p>(1-\log{\pi}/\log{4})^{-1}$ there exist multiplicative Hankel forms in the Schatten class $\mathcal{S}_p$ which lack bounded symbols. The lower bound on $p$ is in a certain sense optimal when the symbol of the multiplicative Hankel form is a product of homogeneous linear polynomials.
\end{abstract}

\maketitle

% INTRODUCTION
\section{Introduction}
For a sequence $\varrho = (\varrho_1,\,\varrho_2,\,\varrho_3,\,\ldots\,) \in \ell^2$ its corresponding \emph{multiplicative Hankel form} on $\ell^2 \times \ell^2$ is given by
\begin{equation} \label{eq:multhankelform}
	\varrho(a,b) = \sum_{m=1}^\infty \sum_{n=1}^\infty \varrho_{mn} a_m b_n,
\end{equation}
which initially is defined at least for finitely supported $a,b\in\ell^2$. Such forms are naturally understood as small Hankel operators on the Hardy space of the infinite polydisc, $H^2(\mathbb{D}^\infty)$. Therefore, one is led to investigate the relationship between the symbol --- a function on the polytorus $\mathbb{T}^\infty$ generating the Hankel form --- and the properties of the corresponding Hankel operator.

In the classical setting, (additive) Hankel forms are realized as Hankel operators on the Hardy space in the unit disc, $H^2(\mathbb{D})$. Nehari's theorem \cite{Nehari} states that every bounded Hankel form is generated by a bounded symbol on the torus $\mathbb{T}$. 

On the infinite polydisc, the study of the corresponding statement was initiated by H. Helson \cite[pp. 52--54]{helsonbook}, who raised the following questions.

\begin{que}
	Does every bounded multiplicative Hankel form have a bounded symbol $\psi$ on the polytorus $\mathbb{T}^\infty$?
\end{que}

\begin{que}
	 Does every multiplicative Hankel form in the Hilbert--Schmidt class $\mathcal{S}_2$ have a bounded symbol?
\end{que}

Helson himself \cite{helson06} gave a positive answer to Question~2. Ortega-Cerd{\`a} and Seip \cite{ortegaseip} proved that there are bounded multiplicative Hankel forms that do not have bounded symbols, using an idea of Helson \cite{helson10}, and hence gave a negative answer to Question~1. Furthermore, their argument also quickly produces that there are compact Hankel forms without bounded symbols (see Lemma~\ref{lem:lifters}). In light of these results, a next natural question to ask is:

\begin{que}	
	Does there exist a Hankel form belonging to a Schatten class $\mathcal{S}_p$, $2 < p < \infty$, without a bounded symbol? If so, for which values of $p$ does such a form exist?
\end{que}

We will answer the first part of this question, by showing that for every
\[p>p_0 = \left(1-\frac{\log{\pi}}{\log{4}}\right)^{-1}\approx 5.738817179,\]
there are multiplicative Hankel forms in $\mathcal{S}_p$ which do not have bounded symbols. 

Our construction relies on independent products of homogeneous linear symbols and is optimal when testing against products of linear homogeneous polynomials, see Theorem~\ref{thm:optimal}. It is quite tempting to further conjecture that forms without bounded symbols can be found in $\mathcal{S}_p$ for every $p>2$, but our method does not substantiate this claim.

The paper is organized into two further sections. Section~2 reviews the connection between multiplicative Hankel forms, the Hardy space of Dirichlet series, and the Hardy space of the infinite polydisc. In Section~3 the main results are proven. 
% The  corresponding statement on the infinite polydisc was first given attention by Helson \cite{helson06}, and in \cite{helson10} he established that multiplicative Hankel forms in the Hilbert-Schmidt class $\mathcal{S}_2$ do have bounded symbols. However, the full analogue of Nehari's theorem was shown to be false by Ortega-Cerd\`a and Seip \cite{ortegaseip} -- there are plenty of bounded Hankel forms which do not have a bounded symbol.

%We want to address the Hankel forms with a degree of smoothness lying in between bounded and Hilbert-Schmidt. The argument of Ortega-Cerd\`a and Seip quickly produces that there are compact Hankel forms without bounded symbols (see Lemma \ref{lem:lifters}), so the first natural question to ask is about the Schatten classes $\mathcal{S}_p$, $p > 2$. This question is the concern of the present paper. Our main result, Theorem \ref{thm:main}, states that for every $p > p_0$,
%\begin{equation} \label{eq:p0}
%	p_0 = \left(1-\frac{\log{\pi}}{\log{4}}\right)^{-1}\approx 5.738817179,
%\end{equation}
%there is a Hankel form belonging to $\mathcal{S}_p$ which does not have a bounded symbol. 

% PRELIMINARIES
\section{Preliminaries}
We let $\mathscr{H}^2$ denote the Hilbert space of Dirichlet series
\begin{equation} \label{eq:diriseri}
	f(s) = \sum_{n=1}^\infty a_n n^{-s}
\end{equation}
with square summable coefficients. If $g$ and $\varphi$ are Dirichlet series in $\mathscr{H}^2$ with coefficients $b_n$ and $\overline{\varrho_n}$, respectively, a computation shows that
\[ \langle fg, \varphi \rangle_{\mathscr{H}^2} = \varrho(a,b).\] 
A key tool in the study of Hardy spaces of Dirichlet series is the \emph{Bohr lift} \cite{bohr}. For any $n\in\mathbb{N}$, the fundamental theorem of arithmetic yields the prime factorization
\[n = \prod_{j=1}^{\infty} p_j^{\kappa_j},\]
which associates the finite non-negative multi-index $\kappa(n) = (\kappa_1,\,\kappa_2,\,\kappa_3,\,\ldots\,)$ to $n$. The Bohr lift of the Dirichlet series \eqref{eq:diriseri} is the power series
\begin{equation} \label{eq:bohrlift}
	\mathscr{B}f(z) = \sum_{n=1}^\infty a_n z^{\kappa(n)},
\end{equation}
where $z = (z_1,\,z_2,\,z_3,\,\ldots\,)$. Hence \eqref{eq:bohrlift} is a power series in countably infinite number of variables, but each term contains only a finite number of variables. 

Under the Bohr lift, $\mathscr{H}^2$ corresponds to the infinite dimensional Hardy space $H^2(\mathbb{D}^\infty)$, which we view as a subspace of $L^2(\mathbb{T}^\infty)$. We refer to \cite{HLS} for the details, mentioning only that the Haar measure of the compact abelian group $\mathbb{T}^\infty$ is simply the product of the normalized Lebesgue measures of each variable. In particular, $H^2(\mathbb{D}^d)$ is a natural subspace of $H^2(\mathbb{D}^\infty)$. 

A formal computation shows that 
\[\langle \mathscr{B}f\mathscr{B}g,\mathscr{B}\varphi\rangle_{L^2(\mathbb{T}^\infty)} = \langle fg,\varphi\rangle_{\mathscr{H}^2},\]
allowing us to compute the multiplicative Hankel form \eqref{eq:multhankelform} on $\mathbb{T}^\infty$. In the remainder of this paper we work exclusively in the polydisc, with no reference to Dirichlet series. Therefore, we drop the notation $\mathscr{B}$ and study Hankel forms
\begin{equation} \label{eq:polyform}
H_\varphi(fg) = \langle f g,\varphi\rangle_{L^2(\mathbb{T}^\infty)}, \qquad f, g \in H^2(\mathbb{D}^\infty).
\end{equation}
In the previous considerations we had that $\varphi \in H^2(\mathbb{D}^\infty)$, but there is nothing to prevent us from considering arbitrary symbols from $L^2(\mathbb{T}^\infty)$. Hence, each $\varphi \in L^2(\mathbb{T}^\infty)$ induces by \eqref{eq:polyform} a (possibly unbounded) Hankel form $H_\varphi$ on $H^2(\mathbb{D}^\infty) \times  H^2(\mathbb{D}^\infty)$. Of course, this is not a real generalization. Each form $H_\varphi$ is also induced by a symbol $\psi \in H^2(\mathbb{D}^\infty)$; letting $\psi = P \varphi$ we have $H_\varphi = H_\psi$, where $P$ denotes the orthogonal projection of $L^2(\mathbb{T}^\infty)$ onto $H^2(\mathbb{D}^\infty)$.
 
Note that if $\psi \in L^\infty(\mathbb{T}^\infty)$, then the corresponding multiplicative Hankel form is bounded, since
\[|H_\psi(fg)| = |\langle fg, \psi \rangle|  \leq \|f\|_2 \, \|g\|_2\|\,\psi\|_\infty.\]
We say that $H_\varphi$ has a bounded symbol if there exists a $\psi \in L^\infty(\mathbb{T}^\infty)$ such that $H_\varphi = H_\psi$. As mentioned in the introduction, it was shown in \cite{ortegaseip} that not every bounded multiplicative Hankel form has a bounded symbol.

On the polydisc the Hankel form $H_\varphi$ is naturally realized as a (small) Hankel operator $\mathbf{H}_\varphi$, which when bounded acts as an operator from $H^2(\mathbb{D}^\infty)$ to the anti-analytic space $\overline{H^2}(\mathbb{D}^\infty)$. Letting $\overline{P}$ denote the orthogonal projection of $L^2(\mathbb{T}^\infty)$ onto $\overline{H^2}(\mathbb{D}^\infty)$, we have at least for polynomials $f \in H^2(\mathbb{D}^\infty)$ that
\begin{equation} \label{eq:operdef}
	\mathbf{H}_\varphi f = \overline{P}( \overline{\varphi} f).
\end{equation}
It is clear that when written in standard bases, the form $H_\varphi$ and the operator $\mathbf{H}_\varphi$ both correspond to the same infinite matrix 
\[M_\varrho = \begin{pmatrix}
			\varrho_{1} & \varrho_{2} & \varrho_{3} & \cdots \\
			\varrho_{2} & \varrho_{4} & \varrho_{6} & \cdots \\
			\varrho_{3} & \varrho_{6} & \varrho_{9} & \cdots \\
			\vdots      & \vdots      & \vdots      & \ddots
		\end{pmatrix}.\]
Finally, we briefly recall the definition of the Schatten classes $\mathcal{S}_p$, $0 < p < \infty$. Assume that the Hankel form $H_\varphi$ is compact. Let $\Lambda = \{\lambda_k\}_{k=1}^\infty$ denote the singular value sequence of $M_\varrho$, which of course is the same as the singular value sequence of the operator $\mathbf{H}_\varphi$. The form $H_\varphi$, or equivalently the operator $\mathbf{H}_\varphi$, is in the \emph{Schatten class} $\mathcal{S}_p$ if $\Lambda \in \ell^p$, and 
\[\|H_\varphi\|_{\mathcal{S}_p}= \|\mathbf{H}_\varphi\|_{\mathcal{S}_p} = \|\Lambda\|_{\ell^p}.\]

% RESULTS  
\section{Results}
To prove that there for each $p > p_0$ exist multiplicative Hankel forms in $\mathcal{S}_p$ without bounded symbols, we will assume that every $H_\varphi \in \mathcal{S}_p$ has a bounded symbol and derive a contradiction. We begin with the following routine lemma.

\begin{lem} \label{lem:lifters}
	Let $p\geq 1$. Assume that every $H_\varphi \in \mathcal{S}_p$ has a bounded symbol on $\mathbb{T}^\infty$. Then there is a constant $C_p\geq1$ with the property that every $H_\varphi \in \mathcal{S}_p$ has a symbol $\psi \in L^\infty(\mathbb{T}^\infty)$ with $H_\varphi = H_\psi$ and such that $\|\psi\|_\infty \leq C_p \|H_\varphi\|_{\mathcal{S}_p}$.
	\begin{proof}
		We will define a lifting operator and show that it has to be continuous by appealing to the closed graph theorem. 
		
		Let $\mathrm{BH}$ denote the space of bounded multiplicative Hankel forms. By a standard argument it is isomorphic to the dual space of the weak product $\mathscr{H}^2 \odot \mathscr{H}^2$ \cite{helson10}. In particular $\mathrm{BH}$ is a Banach space under the operator norm. It follows that $\mathcal{S}_p \mathrm{H}$ is also a Banach space, where $\mathcal{S}_p \mathrm{H}$ denotes the space of multiplicative Hankel forms in $\mathcal{S}_p$ equipped with the norm of $\mathcal{S}_p$.
		
		Now we define
		\begin{align*}
			X &= L^\infty(\mathbb{T}^\infty) \cap \left( L^2(\mathbb{T}^\infty)\ominus H^2(\mathbb{T}^\infty)\right), \\
			Y &= L^\infty(\mathbb{T}^\infty) / X.
		\end{align*}
		$Y$ is a Banach space under the norm $\|\varphi\|_Y = \inf\left\{\|\psi\|_\infty \,:\, \psi - \varphi \in X\right\}$, seeing as $X$ is a closed subspace of $L^\infty(\mathbb{T}^\infty)$. Since by assumption every $H_\varphi \in \mathcal{S}_p\mathrm{H}$ has a symbol $\psi \in L^\infty(\mathbb{T}^\infty)$, we can define a map $T\colon \mathcal{S}_p\mathrm{H} \to Y$ by $T(H_\varphi) = \psi$. This is a well-defined linear map since $H_\varphi = 0$ for a symbol $\varphi \in L^\infty(\mathbb{T}^\infty)$ if and only if $\varphi \in X$. An obvious computation verifies that $T$ is a closed operator, hence continuous. Therefore, there is a $C_p\geq1$ such that
		\[\|T(H_\varphi)\|_Y \leq C_p \|H_\varphi\|_{\mathcal{S}_p}.\]
		The statement of the lemma follows immediately.
	\end{proof}
\end{lem}

Given the assumption of the lemma, we hence have for each polynomial $f$ and form $H_\varphi\in\mathcal{S}_p$ that
\[|\langle f, \varphi \rangle| = |H_\varphi(f\cdot1)| = |H_\psi(f\cdot1)| = |\langle f, \psi \rangle| \leq \|\psi\|_\infty \|f\|_1 \leq C_p \|H_\varphi\|_{\mathcal{S}_p} \|f\|_1,\]
where $\| \cdot \|_1$ denotes the norm of $L^1(\mathbb{T}^\infty)$. We thus obtain
\begin{equation}
	\label{eq:Cest}
	\frac{|\langle f, \varphi \rangle|}{\|H_\varphi\|_{\mathcal{S}_p} \, \|f\|_1} \leq C_p
\end{equation}
for every polynomial $f$ and every $H_\varphi\in\mathcal{S}_p$. To prove our main result we will construct a sequence of polynomials and finite rank forms to show that no finite constant $C_p$ satisfying \eqref{eq:Cest} exists for $p>p_0$, thus obtaining a contradiction to the assumption of Lemma~\ref{lem:lifters}. We will require the following lemma.

\begin{lem} \label{lem:Ssplit}
	Suppose that $\varphi_1,\,\varphi_2,\,\ldots,\,\varphi_m$ are symbols that depend on mutually separate variables and which generate the multiplicative Hankel forms $H_{\varphi_j} \in \mathcal{S}_p$, $1 \leq j \leq m$. Then 
	\begin{equation} \label{eq:sptensor}
	\|H_\varphi\|_{\mathcal{S}_p} = \|H_{\varphi_1}\|_{\mathcal{S}_p}\,\|H_{\varphi_2}\|_{\mathcal{S}_p}\, \cdots \,\|H_{\varphi_m}\|_{\mathcal{S}_p},
	\end{equation}
	where $\varphi = \varphi_1\varphi_2\cdots\varphi_m$.
	\begin{proof}
		For $1 \leq j \leq m$, we let $X_j$ denote the Hardy space of precisely the variables that the symbol $\varphi_j$ depends on, and if necessary let $X_0$ denote the Hardy space of the remaining variables, so that --- as tensor products of Hilbert spaces --- we have
		\[H^2(\mathbb{D}^\infty) = X_0 \otimes X_1 \otimes X_2 \otimes \cdots X_m.\]
		We set $\varphi_0=1$ and consider the small Hankel operators $\widetilde{\mathbf{H}}_{\varphi_j} \colon X_j \to \overline{X_j}$, defined similarly to \eqref{eq:operdef} for $0 \leq j \leq m$. Now, if $f_j \in X_j$, $0 \leq j \leq m$, we observe that
		\[\mathbf{H}_\varphi(f_0f_1\cdots f_m)=\widetilde{\mathbf{H}}_{\varphi_0}(f_0)\,\widetilde{\mathbf{H}}_{\varphi_1}(f_1)\,\cdots\,\widetilde{\mathbf{H}}_{\varphi_m}(f_m),\]
		and hence $\mathbf{H}_\varphi = \widetilde{\mathbf{H}}_{\varphi_0}\otimes \widetilde{\mathbf{H}}_{\varphi_1} \otimes \cdots \otimes \widetilde{\mathbf{H}}_{\varphi_m}$. 
		
		Note that $\widetilde{\mathbf{H}}_{\varphi_0}$ has the sole singular value $1$, of multiplicity $1$. It follows that all singular values $\lambda$ of $\mathbf{H}_{\varphi}$ are obtained as products $\lambda = \lambda_1 \lambda_2 \cdots \lambda_m$, where $\lambda_j$ is a singular value of $\widetilde{\mathbf{H}}_{\varphi_j}$, see \cite{bp}. The multiplicity of $\lambda$ is also obtained in the expected way. From this, a short computation shows that
		\[\|\mathbf{H}_\varphi\|_{\mathcal{S}_p} = \|\widetilde{\mathbf{H}}_{\varphi_1}\|_{\mathcal{S}_p}\,\|\widetilde{\mathbf{H}}_{\varphi_2}\|_{\mathcal{S}_p}\,\cdots\, \|\widetilde{\mathbf{H}}_{\varphi_m}\|_{\mathcal{S}_p}.\]
		Finally, we have $\mathbf{H}_{\varphi_j}=\widetilde{\mathbf{H}}_{\varphi_0}\otimes \widetilde{\mathbf{H}}_{\varphi_j}$, where we now regard $\widetilde{\mathbf{H}}_{\varphi_0}$ as an operator on the Hardy space of the variables of which $\varphi_j$ is independent. Arguing as above, it follows that $\|\mathbf{H}_{\varphi_j}\|_{\mathcal{S}_p}=\|\widetilde{\mathbf{H}}_{\varphi_j}\|_{\mathcal{S}_p}$, completing the proof.
	\end{proof}
\end{lem}

 If $f_1,\,f_2,\,\ldots,\,f_m$ are polynomials depending on the same separate variables as $\varphi_1,\,\varphi_2,\,\ldots,\,\varphi_m$, respectively, and we set $f = f_1f_2\cdots f_m$, then
\begin{align*}
	|\langle f , \varphi \rangle| &= |\langle f_1 , \varphi_1\rangle| \, |\langle f_2 , \varphi_2\rangle|\cdots|\langle f_m , \varphi_m\rangle|, \\
	\|f\|_1 &= \|f_1\|_1\,\|f_2\|_1\cdots \|f_m\|_1.
\end{align*}
Let $\mathrm{S}$ be the shift operator $\mathrm{S}f(z_1, z_2, \ldots) = f(z_2, z_3, \ldots)$. Suppose that we can find polynomials $f$ and $\varphi$, both depending on the first $d$ variables $z_1, z_2, \ldots, z_d$, satisfying
\begin{equation} \label{eq:whatwewant}
	\frac{|\langle f, \varphi \rangle|}{\|H_\varphi\|_{\mathcal{S}_p}\, \|f\|_1}>1.
\end{equation}
Then, for $1 \leq j \leq m$, consider the functions
\[\varphi_j(z) = \mathrm{S}^{d(j-1)}\varphi(z) \quad \text{ and }\quad f_j(z) = \mathrm{S}^{d(j-1)}f(z).\]
With $\Phi = \varphi_1 \varphi_2 \cdots \varphi_m$ and $F = f_1 f_2 \cdots f_m$, Lemma~\ref{lem:Ssplit}  yields
\[\frac{|\langle F, \Phi \rangle|}{\|H_{\Phi}\|_{\mathcal{S}_p}\, \|F\|_1} = \left(\frac{|\langle f, \varphi \rangle|}{\|H_\varphi\|_{\mathcal{S}_p}\, \|f\|_1}\right)^m \to \infty, \qquad m \to \infty,\]
giving us the sought contradiction to \eqref{eq:Cest}. We realize this scheme in the next theorem.

\begin{thm} \label{thm:main}
	For every $p>p_0$ there is a multiplicative Hankel form $H_\varphi \in \mathcal{S}_p$ which does not have a bounded symbol.
	\begin{proof}
		Let $d$ be a large positive integer to be chosen later. Consider the symbol
		\[\varphi (z)= \frac{z_1 + z_2 + z_3 + \cdots + z_d}{\sqrt{d}}.\]
		It is clear that the sequence $\varrho = (\varrho_n)_{n=1}^\infty$ for the matrix of $H_\varphi$ is given by
		\[\varrho_{n} = \begin{cases}
			1/\sqrt{d} & \text{if } n = p_j \text{ and } 1 \leq j \leq d \\
			0 & \text{otherwise}
		\end{cases},\]
		where $p_j$ denotes the $j$th prime. In other terms, the matrix $M_\varrho$ of $H_\varphi$, with all zero rows and columns omitted, is the $(d+1)\times(d+1)$ matrix
		\[\frac{1}{\sqrt{d}}\begin{pmatrix}
			0 & 1 & 1 & \cdots & 1 \\
			1 & 0 & 0 & \cdots & 0 \\
			1 & 0 & 0 & \cdots & 0 \\
			\vdots & \vdots & \vdots & \ddots & \vdots \\
			1 & 0 & 0 & \cdots & 0
		\end{pmatrix}.\]
		This matrix is easily seen to have the singular values $1$ (with multiplicity $2$) and $0$ (with multiplicity $d-1$), and thus
		\[\|H_{\varphi}\|_{\mathcal{S}_p} = 2^{\frac{1}{p}}.\]		
		We choose $f(z) = \varphi(z)$. Then $\langle f,\varphi\rangle=1$, and, moreover, the central limit theorem for Steinhaus variables gives us that
		\[\lim_{d \to \infty} \|f\|_1 = \lim_{d \to \infty} \mathbb{E} \left(\frac{|z_1 + z_2 + z_3 + \cdots + z_d|}{\sqrt{d}}\right) = \frac{\sqrt{\pi}}{2}.\]
		In particular, for each $\delta>0$ we have for sufficiently large $d$  that
		\[\|f\|_1 \leq \frac{\sqrt{\pi}}{2}+\delta.\]
		We now observe that $p=p_0$ is the solution of the equation $2^{1/p} \cdot \sqrt{\pi}/2=1$, and hence if $p>p_0$ we may find $\delta>0$ small enough that
		\[\|H_\varphi\|_{\mathcal{S}_p}\cdot \|f\|_1 \leq 2^{1/p} \cdot\left( \frac{\sqrt{\pi}}{2}+\delta\right)<1.\]
		This implies that if $d$ is large enough, $f$ and $\varphi$ satisfy \eqref{eq:whatwewant}. This completes the proof by appealing to the discussion preceding the statement of the theorem.
	\end{proof}
\end{thm}

Our result is optimal for symbols which are independent products of linear homogeneous polynomials and test functions of the same form, as shown by the following result.
\begin{thm} \label{thm:optimal}
	Suppose $p\leq p_0$ and consider
	\begin{align*}
		\varphi(z) &= a_1 z_1 + a_2 z_2 + \cdots + a_d z_d \quad\text{ and }\quad f(z) = b_1 z_1 + b_2 z_2 + \cdots + b_d z_d,
	\end{align*}
	for $a_j,\,b_j \in \mathbb{C}$. Then $|\langle f, \varphi \rangle | \leq \|H_\varphi\|_{\mathcal{S}_p} \|f\|_1$. 
	\begin{proof}
		By the Cauchy--Schwarz inequality and Parseval's formula, it is clear that
		\[\left|\langle f, \varphi \rangle \right| \leq \|a\|_{\ell^2} \|b \|_{\ell^2}.\]
		 Straightforward computations with the matrix $M_\varrho$ of $H_\varphi$ show that
		\[M_\varrho M_\varrho^\ast = \begin{pmatrix}
			\|a\|_{\ell^2}^2 & 0 & 0 & \cdots & 0 \\
			0 & a_1\overline{a_1} & a_1 \overline{a_2} & \cdots & a_1 \overline{a_d} \\
			0 & a_2\overline{a_1} & a_2 \overline{a_2} & \cdots & a_2 \overline{a_d} \\
			\vdots & \vdots & \vdots & \ddots & \vdots \\
			0 & a_d\overline{a_1} & a_d \overline{a_2} & \cdots & a_d \overline{a_d} \\
		\end{pmatrix}.\]
		Here we have again omitted zero rows and columns. Note that the lower right block has rank $1$. By considering the vector $(0,\,a_1,\,a_2,\,\ldots,\,a_d)$ it is clear that it has the sole eigenvalue $\|a\|_{\ell^2}^2$. Thus, the singular value sequence of $M_\varrho$ is $\Lambda = \{\|a\|_{\ell^2},\,\|a\|_{\ell^2},\,0,\,\ldots,\,0\}$, and hence
		\[\|H_\varphi\|_{\mathcal{S}_p}=2^{1/p}\|a\|_{\ell^2}.\]
		We use the optimal Khintchine inequality for Steinhaus variables \cite{konig,sawa}, $p=1$, and obtain
		\[\|f\|_1 \geq \frac{\sqrt{\pi}}{2} \|b\|_{\ell^2}.\]
		The hypothesis that $p\leq p_0$ implies that $2^{1/p}\sqrt{\pi}/2 \geq 1$, and the proof is finished by the following chain of inequalities.
		\[\|H_\varphi\|_{\mathcal{S}_p}\cdot \|f\|_1 \geq 2^{1/p}\cdot\|a\|_{\ell^2}\cdot\frac{\sqrt{\pi}}{2}\cdot\|b\|_{\ell^2} \geq \|a\|_{\ell^2}\cdot \|b\|_{\ell^2} \geq |\langle f, \varphi\rangle |. \qedhere\]
	\end{proof}
\end{thm}

% REFERENCES
\bibliographystyle{amsplain} 
\bibliography{hankel2} 

\providecommand{\bysame}{\leavevmode\hbox to3em{\hrulefill}\thinspace}
\providecommand{\MR}{\relax\ifhmode\unskip\space\fi MR }
\providecommand{\MRhref}[2]{%
  \href{http://www.ams.org/mathscinet-getitem?mr=#1}{#2}
}
\providecommand{\href}[2]{#2}
\begin{thebibliography}{10}

\bibitem{bohr}
H.~Bohr, \emph{Ueber die {B}edeutung der {P}otenzreihen unendlich vieler
  {V}ariablen in der {T}heorie der {D}irichletschen {R}eihe}, Nachr. Akad.
  Wiss. G{\"o}ttingen Math.--Phys. \textbf{1913} (1913), 441--488.

\bibitem{bp}
A.~Brown and C.~Pearcy, \emph{Spectra of tensor products of operators}, Proc.
  Amer. Math. Soc. \textbf{17} (1966), no.~1, 162--166.

\bibitem{HLS}
H.~Hedenmalm, P.~Lindqvist, and K.~Seip, \emph{A {H}ilbert space of {D}irichlet
  series and systems of dilated functions in {$L^2(0,1)$}}, Duke Math. J.
  \textbf{86} (1997), no.~1, 1--37.

\bibitem{helsonbook}
H.~Helson, \emph{Dirichlet series}, Henry Helson, 2005.

\bibitem{helson06}
\bysame, \emph{Hankel forms and sums of random variables}, Studia Math.
  \textbf{176} (2006), 85--92.

\bibitem{helson10}
\bysame, \emph{Hankel forms}, Studia Math. \textbf{198} (2010), 79--84.

\bibitem{konig}
H.~K{\"o}nig, \emph{On the best constants in the {K}hintchine inequality for
  {S}teinhaus variables}, Israel J. Math. (2013), 1--35.

\bibitem{Nehari}
Z.~Nehari, \emph{On bounded bilinear forms}, Ann. of Math. (1957), 153--162.

\bibitem{ortegaseip}
J.~Ortega-Cerd{\`a} and K.~Seip, \emph{A lower bound in {Nehari’s} theorem on
  the polydisc}, J. Anal. Math. \textbf{118} (2012), no.~1, 339--342.

\bibitem{sawa}
J.~Sawa, \emph{The best constant in the {K}hintchine inequality for complex
  {S}teinhaus variables, the case $p=1$}, Studia Math. \textbf{81} (1985),
  no.~1, 107--126.

\end{thebibliography}
\end{document}